\documentclass[12pt,a4paper, twoside]{article}

\usepackage{amsmath,amssymb,amsthm,amsfonts,mathrsfs,amscd,environ}
\usepackage{dsfont}
\usepackage{latexsym,enumerate,color,geometry,extarrows}
\usepackage{xcolor}
\usepackage{verbatim,fancyhdr,enumitem}
 \NewEnviron{ews}{%
\begin{equation}\begin{split}
  \BODY
\end{split}\end{equation}
}

\NewEnviron{ews*}{%
\begin{equation*}\begin{split}
  \BODY
\end{split}\end{equation*}
}

\def\beg{\begin}
\def\bequ{\begin{equation}}
\def\enqu{\end{equation}}
\def\bes{\begin{split}}
\def\ens{\end{split}}
\def\bews{\begin{ews}}
\def\beqn{\begin{eqnarray}}
\def\enqn{\end{eqnarray}}
\def\beq*{\begin{equation*}}
\def\enq*{\end{equation*}}
\def\bqn*{\begin{eqnarray*}}
\def\eqn*{\end{eqnarray*}}
\def\bary{\begin{array}}
\def\eary{\end{array}}
\def\bpma{\begin{pmatrix}}
\def\epma{\end{pmatrix}}
\def\bvma{\begin{Vmatrix}}
\def\evma{\end{Vmatrix}}

 \numberwithin{equation}{section}

\def\al{\alpha}
\def\be{\beta}

\def\de{\delta}
\def\ep{\epsilon}

\def\et{\eta}

\def\ka{\kappa}
\def\la{\lambda}
\def\rh{\rho}

\def\si{\sigma}
\def\ta{\tau}
\def\ph{\phi}

\def\ps{\psi}

\def\tri{\triangle}
\def\til{\tilde}

\def\Ga{\Gamma}

\def\De{\Delta}

\def\Ph{\Phi}

\def\Om{\Omega}

\def\Q{\mathbb Q}
\def\R{\mathbb R}
\def\P{\mathbb P}
\def\E{\mathbb E}

\def\N{\mathbb N}

\def\W{\mathbb W}

\def\sF{\mathscr F}
\def\sD{\mathscr D}
\def\sC{\mathscr C}
\def\sB{\mathscr B}
\def\sH{\mathscr H}
\def\sL{\mathscr L}

\def\sP{\mathscr P}

\def\cF{\mathcal F}

\def\e{\operatorname{e}}

\def\d{\mathrm{d}}

\def\ff{\frac}
\def\ra{\rightarrow}
\def\da{\downarrow}
\def\na{\nabla}
\def\pp{\partial}
\def\<{\langle}
\def\>{\rangle}
\def\sq{\sqrt}
\def\tld{\tilde}
\def\we{\wedge}
\def\1{\mathds{1}}

\def\trac{\mathrm{tr}}

\def\8{\infty}

\allowdisplaybreaks

\title{{\bf TCI for SDEs with irregular drifts}
}

\author{
{\bf Yongqiang Suo$^{a)}$,~Chenggui Yuan$^{a)}$~and~Shao-Qin Zhang$^{b)}$}\\
~\\
\footnotesize{$^{a)}$Department of Mathematics,   Swansea University, Bay Campus SA1 8EN, UK.}\\
\footnotesize{ Emails: 971001@swansea.ac.uk$\qquad$C.Yuan@swansea.ac.uk}\\
\footnotesize{$^{b)}$School of Statistics and Mathematics}\\
\footnotesize{Central University of Finance and Economics, Beijing 100081, China}\\
\footnotesize{Email: zhangsq@cufe.edu.cn}\\
}

\begin{document}

\maketitle

\begin{abstract}
We obtain $T_2(C)$ for  stochastic differential equations with Dini continuous drift and $T_1(C)$ stochastic differential equations  with singular coefficients.
\end{abstract}\noindent

\noindent


\vskip 2cm

\section{Introduction}
Let $(E,\rho)$ be a metric space equipped with a $\si$-field $\sB$ such that $\rho(\cdot,\cdot)$ is $\sB\times\sB$ measurable and let $\sP(E)$ be the class of all probability measures on $E$. The $p$-th Wasserstein distance between $\mu,\nu\in\sP(E)$ is defined by 
\beg{align*}
\W_p^\rho(\mu,\nu)=\inf_{\pi\in\sC(\mu,\nu)}\Big(\int_{E\times E}\rho^p(x,y)\pi(\d x,\d y)\Big)^{\ff{1}{p\vee 1}},
\end{align*}
where $\sC(\mu,\nu)$ is the space of all couplings of $\mu,\nu$.
The relative entropy of $\mu$ with respect to $\nu$ is given by
\beg{align*}
H(\nu|\mu)=\begin{cases}\int_E\ln\ff{\d\nu}{\d\mu}\d\nu,\quad \ \ & 
 \mbox{if}~~\nu<<\mu,\\ 
  +\8, \quad \ \ & \mbox{otherwise}.
\end{cases} 
\end{align*}
We say that the probability measure $\mu$ satisfies the $W_p$-transportation cost-information inequality ( TCI for short) on $(E,\rho)$ if there exists a constant $C>0$ such that for any probability measure $\nu$, 
\beg{align*}
W_p^\rho(\mu,\nu)\le \sq{2CH(\nu|\mu)}.
\end{align*}
To  be short, we write $\mu\in T_p(C)$ for this relation.

Since Talagrand's work \cite{T}, the $T_1(C)$   and the $T_2(C)$ have been intensively investigated and applied to many other distributions, such as \cite{DGW,WZ1,WZ2} for diffusion processes, \cite{M,S,W} for stochastic differential equations (SDEs) with L\'evy noise or fractional Brownian motion, \cite{BWY,U} for stochastic functional differential equations (SFDEs). The $T_2(C)$ intimately linked to the concentration of measure phenomenon and some functional inequalities  such as Poincar\'e inequality, log-Sobolev inequality and Hamiton-Jacobi equations, see \cite{FS,GRS,L,OV,WF,WZ2} and references therein. For example, $T_2(C)$ can be derived from the log-Sobolev inequality \cite{BGL,OV}.  Additionally, $T_2(C)$ implies the Poincar\'e inequality. Moreover, the $T_2(C)$ can also be established when the log-Sobolev inequality is unknown, see for instance \cite{BWY,DGW} and references therein. $T_2(C)$ implies the concentration of measure, the converse implication also holds, i.e. If $\mu$ satisfies the normal concentration, then $\mu$ satisfies $T_2(C),$ see \cite{L,SZ}. As for $T_1(C),$ we highlight that \cite{DGW} gave a characterization of $T_1(C)$ by ``Gaussian tail'' on a metric space and some applications to random dynamic systems and diffusions.  Using Malliavin calculus, \cite{M,W}  proved $T_1(C)$ for invariant probability measure of solution to SDEs with the $L^1$- metric and  uniform metric under dissipative conditions.

It is worth noting that most of the above references of TCIs for solutions to SDEs and SFDEs are required to meet Lipschitz condition for the drifts, some references relaxed this condition to the case with one-sided Lipschitz condition. Motivated by \cite{WFY,WZ}, the goal of this paper is to establish the equivalent expressions of Wasserstein distance and relative entropy of measures defined on a polish space by introducing a Homeomorphism on it, which implies the equivalent expression of $T_p(C)$ for laws of solutions to two equivalent SDEs, the coefficients of one  SDE are
singular.  

The remainder of the paper is organised as follows: In section 2, we present a general result on $T_p(C)$ for measure $\mu$ on Polish space $(E,\rho)$. In Section 3, the main results including the $T_2(C)$ for SDEs with Dini continuous drift and $T_1(C)$ for SDEs with singular  coefficients are introduced. By the general results in Section 2,  $T_2(C)$ for SDE \eqref{eq1}  and  the $T_1(C)$ for SDE \eqref{b-S} are  proved in Section  4 and Section 5, respectively.

\section{A general result}
Let $(E,\rh)$ be a Polish space and $\Ph$ be a homeomorphism on $E$ with positive constants $c_1$ and $c_2$ such that 
\beg{align}\label{inq-1}
c_1\rh(x,y)\leq \rh(\Ph(x),\Ph(y))\leq c_2\rh(x,y),~x,y\in E.
\end{align} 
We can see that $\Ph$ induces a homeomorphism on $E\times E$, which is still denoted by $\Ph$: 
$$\Ph(x,y)=\left(\Ph(x),\Ph(y)\right),~(x,y)\in E\times E.$$
It is clear that the inverse of $\Ph$ on $E\times E$ is given by 
$$\Ph^{-1}(x,y)=\left(\Ph^{-1}(x),\Ph^{-1}(y)\right),~x,y\in E.$$
We can now formulate the following result, which is a simple extension of conclusion of \cite[Lemma 2.1]{DGW}. Here,  we give a detailed proof for readers' convenience.
\beg{lem}\label{gthm}
 For any $p\geq 1, \mu,\nu\in\sP(E)$, we have the following assertions hold.
\beg{enumerate}
\item[(1)] 
\beg{align}\label{equ-1}
\W^\rh_p(\mu,\nu)=\W_p^{\rh\circ \Ph^{-1}}(\mu\circ \Ph^{-1},\nu\circ\Ph^{-1}).
\end{align}
If \eqref{inq-1} holds for metric $\rh\circ\Ph^{-1}$, then 
\beg{align}\label{inq-2}
c_1\W_p^{\rh}(\mu,\nu)\leq \W_p^{\rh}(\mu\circ\Ph^{-1},\nu\circ\Ph^{-1})\leq c_2\W_p^{\rh}(\mu,\nu),
\end{align}
holds for some constants $c_1$ and $c_2$.
\item[(2)] If $\mu\ll\nu$, we have 
\beg{align*}
H(\mu|\nu)=H(\mu\circ \Ph^{-1}|\nu\circ\Ph^{-1}).
\end{align*}
\end{enumerate}
\end{lem}
\beg{proof}

1).  Let $\pi\in\sC(\mu,\nu)$ and $A\in \sB(E)$. Then one has
\beg{align*}
\pi\circ\Ph^{-1}(A\times E)&=\pi\left(\Ph^{-1}(A\times E)\right)=\pi\left(\Ph^{-1}(A)\times\Ph^{-1}(E)\right)\\
&=\pi\left(\Ph^{-1}(A)\times E\right)=\mu(\Ph^{-1}(A))\\
&=\mu\circ\Ph^{-1}(A).
\end{align*}
Similarly, it is easy to see that
\beg{align*}
\pi\circ\Ph^{-1}(E\times A)=\nu\circ\Ph^{-1}(A).
\end{align*}
Thus $\pi\circ\Ph^{-1}\in\sC(\mu\circ\Ph^{-1},\nu\circ\Ph^{-1})$.   

On the other hand, for any $\tld \pi\in \sC(\mu\circ\Ph^{-1},\nu\circ\Ph^{-1})$, we similarly have $\tld \pi\circ\Ph\in\sC(\mu,\nu)$. Moreover, $(\pi\circ\Ph^{-1})\circ\Ph=\pi$. Define 
\beg{align*}
\left(\Ph^{-1}\right)^\#: \pi\ra \pi\circ \Ph^{-1},~\pi\in \sC(\mu,\nu),
\end{align*}
then $(\Ph^{-1})^{\#}$ is a bijection from $\sC(\mu,\nu)$ to $\sC(\mu\circ\Ph^{-1},\nu\circ\Ph^{-1})$ with inverse $\Ph^{\#}$.  

For any $\pi\in \sC(\mu,\nu)$, the bijection $(\Ph^{-1})^{\#}$ implies that
\beg{align}\label{rho}
\W_p^{\rh}(\mu,\nu)^p&\leq \int_{E\times E}\rh^p(x,y)\pi(\d x,\d y)\\\nonumber
&=\int_{E\times E}\rh^p\circ\Ph^{-1}(x,y)\pi\circ \Ph^{-1}(\d x,\d y)\\\nonumber
&=\int_{E\times E}\rh^p\circ\Ph^{-1}(x,y)\left((\Ph^{-1})^{\#}\pi\right)(\d x,\d y),\nonumber
\end{align}
which implies that 
\beg{align*}
\W_p^{\rh}(\mu,\nu)^p&\leq \inf_{\pi\in\sC(\mu,\nu)}\int_{E\times E}\rh\circ\Ph^{-1}(x,y)\left((\Ph^{-1})^{\#}\pi\right)(\d x,\d y)\\
&=\inf_{\tld \pi\in\sC(\mu\circ\Ph^{-1},\nu\circ\Ph^{-1})}\int_{E\times E}\rh^p\circ\Ph^{-1}(x,y)\tld \pi(\d x,\d y)\\
&=\W_p^{\rh\circ\Ph^{-1}}(\mu\circ\Ph^{-1},\nu\circ\Ph^{-1})^p.
\end{align*}
Since $\Ph^{\#}$ is the inverse of $(\Ph^{-1})^{\#}$, we have 
\beg{align*}
\W_p^{\rh\circ\Ph^{-1}}(\mu\circ\Ph^{-1},\nu\circ\Ph^{-1})^p\leq \W_p^{\rh}(\mu,\nu)^p.
\end{align*}
This, together with \eqref{rho}, yields \eqref{equ-1}.

Applying \eqref{inq-1} to $\rho\circ\Ph^{-1}$, one obtains from the definition of $L^p$-Wasserstein distance that
\beg{align*}
c_1\W_p^{\rh\circ\Ph^{-1}}(\mu\circ\Ph^{-1},\nu\circ\Ph^{-1})
&\leq \W_p^{\rh}(\mu\circ\Ph^{-1},\nu\circ\Ph^{-1})\\
&\leq c_2 \W_p^{\rh\circ\Ph^{-1}}(\mu\circ\Ph^{-1},\nu\circ\Ph^{-1}).
\end{align*}
Combining this with \eqref{equ-1}, we obtain \eqref{inq-2}.

2).  We first assume $\mu\ll\nu$.  For any $A\in E$, if
 $\nu\circ\Ph^{-1}(A)=0$, i.e. $\nu(\Ph^{-1}(A))=0$, then one has  
\beg{align*}
\mu\circ\Ph^{-1}(A)=\mu(\Ph^{-1}(A))=0,
\end{align*}
which implies $\mu\circ\Ph^{-1}\ll\nu\circ\Ph^{-1}$. 

Similarly, if $\mu\circ\Ph^{-1}\ll \nu\circ\Ph^{-1}$, then  $\mu\ll\nu$. 

By the definition of pushforward measure, one obtains that for any $\ps\in\sB(E)$
\beg{align*}
\int_E\ps\ff {\d {\mu\circ\Ph^{-1}}} {\d \nu\circ\Ph^{-1}}\d\nu\circ\Ph^{-1}&=\int_E\ps \d\mu\circ\Ph^{-1}=\int_E\ps\circ\Ph\d \mu\\
&=\int_E\ps\circ\Ph\ff{\d\mu}{\d\nu}\d\nu=\int_E\ps(\ff {\d \mu} {\d \nu}\d \nu)\circ\Ph^{-1}\\
&=\int_E\ps(\ff {\d \mu} {\d \nu}\circ\Ph^{-1})\d \nu\circ\Ph^{-1},
\end{align*}
which yields $\ff {\d {\mu\circ\Ph^{-1}}} {\d \nu\circ\Ph^{-1}}=\ff {\d \mu} {\d \nu}\circ\Ph^{-1}$, $\nu\circ\Ph^{-1}$-a.s.. We then can see that
\beg{align*}
H(\mu|\nu)&=\int_E\left(\log\ff {\d \mu} {\d \nu}\right)\d \mu=\int_E \ff {\d \mu} {\d \nu}\left(\log\ff {\d \mu} {\d \nu}\right) \d \nu\\
&=\int_E\left(\left[\ff {\d \mu} {\d \nu}\left(\log\ff {\d \mu} {\d \nu}\right)\right]\circ\Ph^{-1}\right)\d\nu\circ\Ph^{-1}\\
&=\int_E \left[\ff {\d \mu\circ\Ph^{-1}} {\d \nu\circ\Ph^{-1}}\left(\log\ff {\d \mu\circ\Ph^{-1}} {\d \nu\circ\Ph^{-1}}\right)\right]\d\nu\circ\Ph^{-1}\\
&=H(\mu\circ\Ph^{-1}|\nu\circ\Ph^{-1}).
\end{align*}
\end{proof}

Throughout this work, the following notation will be used. $(\R^d,\<\cdot,\cdot\>,|\cdot|)$ denotes the $d$-dimensional Euclidean space, $\R^d\otimes\R^d$ is the family of all $d\times d$ matrices.  For a vector or matrix $v,$ $v^*$ denotes its transpose. Let $\|\cdot\|$ denote the usual operator norm. Fix $T>0$ and set $\|f\|_{T,\8}:=\sup_{t\in[0,T],x\in\R^d}\|f(t,x)\|$ for an operator or vector valued map $f$ on $[0,T]\times\R^d$, $C(\R^d;\R^d)$ means the set of all continuous functions $f:\R^d\rightarrow\R^d$.  Let $C^2(\R^d;\R^d\otimes\R^d)$ be the family of all continuously twice differentiable functions $f:\R^d\rightarrow\R^d\otimes\R^d$. $\na^i,i\in\N$ means the $i$-th order gradient operator. Let  $W_t$ be  a $d$-dimensional Brownian motion defined on   a complete filtration probability space $(\Om,(\sF_t)_{t\ge0},\sF,\P).$

\section{TCI for SDEs with singular coefficients}
In this section, we will first present $T_2(C)$ for SDEs with Dini continuous drift,   then formulate $T_1(C)$ for SDEs with singular dissipative coefficients.
\subsection{$T_2(C)$ for SDEs with Dini continuous drift}
Consider the following SDE with Dini continuous drift
\beg{align}\label{eq1}
\d X_t=\{B_t(X_t)+b_t(X_t)\}\d t+\si_t(X_t)\d W_t,
\end{align}
where  $B,b:[0,T]\times\R^d\rightarrow\R^d$ are measurable, and $\si:[0,T]\times\R^d\rightarrow\R^d\otimes\R^d$ is measurable.

Let
\beg{align*}
\sD=\Big\{\phi:[0,\8)\rightarrow[0,\8)~\mbox{is increasing},~ \phi^2~\mbox{is concave}, \int_0^1\ff{\phi(s)}{s}\d s<\8\Big\}.
\end{align*}
With regard to \eqref{eq1}, we impose the following conditions on its coefficients. For any fixed $T>0, t\in[0,T], x,y\in\R^d$, there exists $\phi\in\sD$ such that
\beg{enumerate}
\item[(A1)] $\|b\|_{T,\8}<+\8$ and there exists $\phi\in\sD$ such that
\beg{align*}
|b_t(x)-b_t(y)|\le\phi(|x-y|),~~t\in[0,T], x,y\in\R^d.
\end{align*}
 \item[(A2)] $B_t(\cdot)$ satisfies Lipschitz condition and  $\sup_{t\in [0,T]}|B_t(0)|<\infty$; $\si_t(x)$ is invertible and $\si_t\in C^2(\R^d;\R^d\otimes\R^d)$ with $\sup_{t\in [0,T]}\|\si_t(0)\|<\infty$;  there exists  some positive increasing function $K\in C([0,\8);(0,\8))$ such that
 \beg{align*}
\| \na B\|_{T,\infty}+\|\si\|_{T,\8}
+\|\na\si\|_{T,\8}+\|\na^2\si\|_{T,\8}+\|(\si\si^*)^{-1}\|_{T,\8}\le K(T).
 \end{align*}
\end{enumerate}
\beg{rem}\label{thm1}
According to \cite[Theorem 1]{WFY},  for any $T>0$, the equation \eqref{eq1} has a  unique  strong solution $(X_t)_{t\in[0,T]}$ under the  assumptions (A1)-(A2).
\end{rem}
\beg{rem}
 The condition $\int_0^1\ff{\phi(s)}{s}\d s<\8$ is well known as Dini condition. If $\phi$ is H\"older continuous with exponent $\al$, then $\phi$ is Dini continuous. In fact, if $\phi(0)=0$ and
 $|\phi(s)-\phi(t)|\le L|s-t|^\al$, then $\int_0^1\ff{\phi(s)}{s}\d s\le\ff{L}{\al}$ holds.
However, there are numerous Dini continuous functions which are not H\"older continuous for any $\al>0$. For instance, 
\beg{align*}
\phi(s)=\begin{cases}
(\log s)^{-2},& s\in(0,\e^{-3}),\\
0,& s=0.
\end{cases}
\end{align*}
 It is easy to check that $\lim_{s\rightarrow0^+}\ff{\phi(s)}{s^\al}=\8$ for any $\al>0$, so $\phi$ is not H\"older continuous, but $\phi$ is Dini continuous. Indeed, $\phi'(s)>0,s\in(0,\e^{-3})$ and $\int_0^1\ff{\phi(s)}{s}\d s<\8$, which implies that $\phi$ is Dini continuous.
 \end{rem}
%
\beg{thm}\label{Dthm}
Suppose the assumptions (A1)-(A2) hold. \\
(1) Let $\P^x$  be the law of the solution $X(\cdot)$  to \eqref{eq1} with initial value $X(0)=x$. Then   the quadratic transportation cost inequality on the space $C([0,T];\R^d)$, i.e.
\beg{align*}
\W_{2}^{\rho_T}(\Q,\P^x)^2\le CH(\Q|\P^x), \Q\in\sP(C([0,T];\R^d))
\end{align*}
holds for some constant $C>0$, where $\rho_T$ denotes the uniform metric on the space $C([0,T];\R^d)$.\\
(2) Let $\mu\in\sP(\R^d)$ and $\P^\mu$  be the law of the solution $X(\cdot)$ to \eqref{eq1} with initial distribution $\mu$ . Then 
 \beg{align}\label{eqk1}
 \W_{2}^{\rho_T}(\Q,\P^\mu)^2\le C_1H(\Q|\P^\mu),~~\Q\in\sP(C([0,T];\R^d))
 \end{align}
 holds for some constant $C_1>0$ if and only if 
 \beg{align}\label{eqk2}
 \W_2^{\rho}(\nu,\mu)^2\le C_2H(\nu|\mu),~~\nu\in\sP(\R^d)
 \end{align}
 holds for some constant $C_2>0$.
\end{thm}
\beg{rem}
Based on \cite[Theorem 2.4]{L}, the conclusion of this theorem implies that $\P^x$ satisfies the concentration property  with
\beg{align*}
\al(r)=\e^{-\ff{1}{C}(r-r_0)^2},~r\ge0,
\end{align*}
where  $r_0=\sq{C\log(2)}$ and  the constant $C$ is same as in the above theorem.
\end{rem}

\subsection{$T_1(C)$ for SDEs with singular dissipative coefficients} 
In this subsection, we consider the following singular SDE
\beg{align}\label{b-S}
\d X_t=b(X_t)\d t+\si(X_t)\d W_t,
\end{align}
where   $b:\R^d\ra\R^d, \si:\R^d\ra\R^d\otimes\R^d$ are Borel measurable functions.
Assume that  the coefficients $b$ and $\si$ satisfy  the following conditions:  
\beg{enumerate}
\item[($H^b$)] Assume $b=b_1+b_2$ such that  $b_1\in L^p(\R^d)$ for some $p>d$,  and one of the following conditions holds for $b_2$\\  
(1) for some $\ka_1,\ka_2,\ka_3>0,$  $r>-1$
\beg{align}\label{diss}
\<x,b_2(x)\>\le-\ka_1|x|^{2+r}+\ka_2, ~~\mbox{and}~~|b_2(x)|\le \ka_3(1+|x|^{1+r});
\end{align}
(2) for some $\ka_4\geq 0$
\beg{align}\label{lin}
|b_2(x)|\leq \ka_4(1+|x|).
\end{align}
\item[($H^\si$)] $\|\na\si\|\in L^p(\R^d)$ with the same $p$ in ($H^b$). There are constants $c_0\ge1$ and $\be\in(0,1)$ such that
\beg{align*}
&c_0^{-1}|\xi|^2\le|\si^*(x)\xi|^2\le c_0|\xi|^2,~~\forall \xi\in\R^d,\\
&\|\si(x)-\si(y)\|\le c_0|x-y|^\be, ~~x,y\in\R^d.
\end{align*}
\end{enumerate}
$b_1$ is called  the singular part and $b_2$ is locally bounded.  According to \cite[Theorem 2.10]{XZ}, under $(H^b)$ and $(H^\si)$, \eqref{b-S} admits a unique strong solution.  We now state the $T_1(C)$ for law of solution to SDE \eqref{b-S} with initial point $x$.
\beg{thm}\label{TCI-1}
Assume assumptions $(H^b)$ and $(H^\si)$ hold. Then the law $\P^x$ of solution $X_t(x)$ to SDE \eqref{b-S} satisfies the $T_1(C)$ for every $x\in\R^d$ on the space $C([0,T];\R^d)$ equipped with the uniform norm $\rh_T$.
\end{thm}

\section{Proof of Theorem \ref{Dthm}}
 \subsection{Regularization representation of the solution to \eqref{eq1}}
By Lemma \ref{gthm}, we establish the $T_2(C)$ for $\P^x$ by constructing a differeomorphism on $C([0,T],\R^d)$. To this end, we construct a transform $\Phi:[0,T]\times\R^d\rightarrow\R^d$ in the spirit of \cite{WFY}.  

In the sequel, we briefly explain how to construct the transform $\Phi$. 

We first decompose $B_t$ into a smooth term and a bounded Lipschitz term. 
\beg{lem}
There exist $\bar B_t\in C^2(\R^d)$ and $\hat B_t$ which is Lipschitz such that $B_t=\bar B_t+\hat B_t$ and 
\beg{align}
\|\na \bar B\|_{T,\infty}+\|\na^{2} \bar B\|_{T,\infty}&<\infty,\nonumber\\
\|\na \bar B\|_{T,\infty}\vee\|\na \hat B\|_{T,\infty}\vee \|\hat B\|_{T,\infty}&\leq \|\na B\|_{T,\infty}.\label{nn-b-B}
\end{align}
\end{lem}
\beg{proof}
Let $r(x)$ be a smooth function supported in $\{x\in\R^d~|~|x|\leq 1\}$ and $\int_{\R^d}r(x)\d x=1$. Set
\beg{align*}
\bar B_t(x)=B_t*r(x),\qquad \hat B_t(x)=B_t(x)-\bar B_t(x).
\end{align*}
Then the assertions of this lemma hold. 

\end{proof} 
Let $\hat b_t(x)=\ph_t(x)+\hat B_t(x)$. Then
\beg{align*}
|\hat b_t(x)-\hat b_t(y)|&\leq \ph(|x-y|)+ \|\na B\|_{T,\infty}(|x-y|\wedge 1)\\
&\leq \ph(|x-y|)+ \|\na B\|_{T,\infty}(|x-y|^{\ff 1 2}\wedge 1)\\
&=:\hat\phi(|x-y|),~x,y\in\R^d.
\end{align*}
Moreover, we have that $\hat \phi\in \sD$. Hence, we use the following assumption instead of (A2) 
\beg{enumerate}
\item[(A2')] $B_t(\cdot)\in C^2(\R^d;\R^d)$ with $\sup_{t\in [0,T]}|B_t(0)|<\infty$; $\si_t(x)$ is invertible and $\si_t\in C^2(\R^d;\R^d\otimes\R^d)$ with $\sup_{t\in [0,T]}\|\si_t(0)\|<\infty$;  there exists  some positive increasing function $K\in C([0,\8);(0,\8))$ such that
 \beg{align*}
\| \na B\|_{T,\infty}+\| \na^2 B\|_{T,\infty}+\|\si\|_{T,\8}
+\|\na\si\|_{T,\8}&\\
+\|\na^2\si\|_{T,\8}+\|(\si\si^*)^{-1}\|_{T,\8}&\le K(T).
 \end{align*}
\end{enumerate}

Consider a backward PDE
\beg{align}\label{PDE1}
\pp_t u_t=-L_t u_t-b_t+\la u_t,~u_T=0,~t\in [0,T],
\end{align}
where $\la>0$  is a parameter and $L_t:=\ff{1}{2}\trac{(\si_t\si_t^*\na^2)}+\na_{B_t}+\na_{b_t}$.  Set $\Phi_t(x)=x+u_t(x)$. Then $\pp_t\Phi_t=B_t(x)-L_t\Phi_t(x)+\la u_t(x)$. By It\^o's formula, we formally have that (see Lemma \ref{lem2} for a proof)
\beg{align}\label{eq2}
\d\Phi_t(X_t) &=\{(\pp_t\Phi_t)(X_t)+L_t\Phi_t(X_t)\}\d t+\na\Phi_t(X_t)\si_t(X_t)\d W_t\nonumber\\
& = (\la u_t(X_t)+B_t(X_t))\d t+\na\Phi_t(X_t)\si_t(X_t)\d W_t,
\end{align}
The irregular term $b_t$ is canceled. $u_t$ can be regular and $\|\na u\|_{T,\infty}<1$ for large enough $\la$, see Lemma \ref{lu}. Then $\Ph_t$ is a differeomorphism.

We investigate \eqref{PDE1} in a weaker form. Let $\{P_{s,t}^0\}_{0\le s\le t}$ is the semigroup associated to the SDE below
\beg{align}\label{eq3}
\d Z_{s,t}^x=B_t(Z_{s,t}^x)\d t+\si_t(Z_{s,t}^x)\d W_t,~~t\ge s, Z_{s,s}^x=x.
\end{align}
It is well known that the equation \eqref{eq3} has a unique solution under assumption (A2'). Then we have
\beg{align*}
P_{s,t}^0f(x)=\E f(Z_{s,t}^x),~~t\ge s\ge 0, x\in\R^d, f\in\sB_b(\R^d).
\end{align*}
The generator of $P_{s,t}^0$ is  $\tilde{L}_t=\ff{1}{2}\trac{(\si_t\si_t^*\na^2)}+\na_{B_t}.$  By using $P_{s,t}^0$, \eqref{PDE1} can be rewritten into the following integral equation 
\beg{align}\label{eq4}
u_s=\int_s^T\e^{-\la(t-s)}P_{s,t}^0\{\na_{b_t}u_t+b_t\}\d t, s\in[0,T].
\end{align}

In the following lemma, we give the gradient estimates for semigroup $P_{s,t}^0$ defined by \eqref{eq3}, which will be used to study the regularity properties of solution $u_s$ to equation \eqref{eq4}. The proof of the following lemma follows from \cite{WFY} completely, and we omit it. 
\beg{lem}\label{lem1}
Fix $T>0$. Assume (A2'). Then the following assertions hold.
\beg{enumerate}
\item[(1)] For any $f\in\sB_b(\R^d)$, $P_{s,t}^0f\in C^2_b({\R^d})$. There exists a positive constant $c$ such that for any $0\le s< t\le T$,
\beg{align}\label{eq10}
|\na P_{s,t}^0f|^2(x)&\le 
\ff{c}{t-s}P_{s,t}^0f^2(x),\\
|\na^2 P_{s,t}^0f|^2(x)&\le 
\ff{c}{(t-s)^2}P_{s,t}^0f^2(x),~x\in\R^d,~f\in \sB_b(\R^d). \label{ine-na2}
\end{align} 
\item[(2)]There exist positive constants $c_1$ and $c_2$ such that for any increasing $\phi:[0,\8)\rightarrow[0,\8)$ with concave $\phi^2$ 
\beg{align}\label{eq2.2}
\|\na^2P_{s,t}^0f\|_\8:=\sup_{x\in\R^d}\|\na^2P_{s,t}^0f(x)\|\le \ff{c_1\phi(c_2\sq{(t-s)})}{t-s},
\end{align}
holds for any $f\in\sB_b(\R^d)$ satisfying
\beg{align*}
|f(x)-f(y)|\le\phi(|x-y|),~~0\le s< t\le T, x,y\in\R^d.
\end{align*}
\end{enumerate}
\end{lem}

The following Lemma focuses on the existence and uniqueness of solution to \eqref{eq4} and gradient estimates of the solution, which is essentially due to \cite[Lemma 2.3]{W}. We include a complete proof for readers' convenience.
\beg{lem}\label{lu}
Assume $\|b\|_{T,\8}<\8$ and (A2'). Let $T>0$ be fixed, then there exists a constant $\la(T)>0$ such that the following assertions hold:
\beg{enumerate}
\item[(1)] For any $\la>\la(T)$, \eqref{eq4} has a unique solution $u\in C([0,T];C_b^1(\R^d;\R^d))$ satisfying 
\beg{align}\label{u1}
\lim_{\la\rightarrow\8}\{\|u\|_{T,\8}+\|\na u\|_{T,\8}\}=0.
\end{align}
\item[(2)] Moreover, if (A1) holds, then we have 
\beg{align}\label{u2}
\lim_{\la\rightarrow\8}\|\na^2u\|_{T,\8}=0.
\end{align} 
\end{enumerate}
\end{lem}
\beg{proof}
(1) Let $\sH=C([0,T];C_b^1(\R^d;\R^d))$, which is a Banach space under the norm $\|u\|_{\sH}:=\|u\|_{T,\8}+\|\na u\|_{T,\8}, u\in\sH$.

For any $u\in\sH$, define the mapping
\beg{align*}
(\Ga u)_s(x)=\int_s^T\e^{-\la(t-s)}P_{s,t}^0\{\na_{b_t(\cdot)}u_t(\cdot)+b_t(\cdot)\}(x)\d t.
\end{align*}
Firstly, we claim that $\Ga \sH\subset\sH$. In fact, for any $u\in\sH$, by  \eqref{eq10}, one has
\beg{align}\label{G1}
\|\Ga u\|_{T,\8}&=\sup_{s\in[0,T],x\in\R^d}\Big|\int_s^T\e^{-\la(t-s)}P_{s,t}^0\{\na_{b_t(\cdot)}u_t(\cdot)+b_t(\cdot)\}(x)\d t\Big|\nonumber\\
&\le\sup_{s\in[0,T]}\Big|\int_s^T\e^{-\la(t-s)}\|b\|_{T,\8}(\|\na u\|_{T,\8}+1)\d t\Big|\nonumber\\
&\le\ff{\|b\|_{T,\8}(\|\na u\|_{T,\8}+1)}{\la}<\8, 
\end{align} 
and
\beg{align}\label{G2}
\|\na\Ga u\|_{T,\8}&=\sup_{s\in[0,T],x\in\R^d}\Big|\int_s^T\e^{-\la(t-s)}\na P_{s,t}^0\{\na_{b_t(\cdot)}u_t(\cdot)+b_t(\cdot)\}(x)\d t\Big|\nonumber\\
&\le c\sup_{s\in[0,T]}\Big|\int_s^T\ff{\e^{-\la(t-s)}}{\sq{t-s}}\|b\|_{T,\8}(\|\na u\|_{T,\8}+1)\d t\nonumber\\
&\le\ff{c\|b\|_{T,\8}(\|\na u\|_{T,\8}+1)}{\sq{\la}}<\8. 
\end{align} 
Therefore,  the claim $\Ga\sH\subset\sH$ holds. 

Next, we will show that for large enough $\la>0$, $\Ga$ is contractive on $\sH$. Indeed, by the similar arguments as above, it is easy to check that for any $u,\hat{u}\in\sH$, we have
\beg{align*}
\|\Ga u-\Ga\hat{u}\|_{\sH}&\le \ff{\|b\|_{T,\8}}{\la}(1+c)\|\na u-\na\hat{u}\|_{T,\8}\\
&\le \ff{\|b\|_{T,\8}}{\la}(1+c)\|u-\hat{u}\|_{\sH}\\
&=:C(\la)\|u-\hat{u}\|_{\sH}.
\end{align*}
Choosing constant $\la(T)$ satisfies $C(\la)<1$ for $\la>\la(T)$, we can see that $\Ga$ is contractive on $\sH$ with $\la>\la(T)$. Thus, the fixed point theorem yields that \eqref{eq4} has a unique solution $u\in\sH$. 

Finally, the estimates \eqref{G1} and \eqref{G2} imply that \eqref{u1} holds.

(2) \eqref{ine-na2} implies that for any $f\in\sB_b(\R^d)$
\beg{align*}
|\na P_{s,t}^0f(x)-\na P_{s,t}^0f(y)|\le\ff{c|x-y|}{t-s}\|f\|_\8,~~x,y\in\R^d,~0\le s< t\le T.
\end{align*}
This, together with \eqref{eq10}, yields that
\beg{align}\label{eq2.1}
|\na P_{s,t}^0f(x)-\na P_{s,t}^0f(y)|\le c\left(\ff{|x-y|}{t-s}\wedge\ff{1}{\sq{t-s}}\right)\|f\|_\8,
\end{align}
where $c$ is some constant.

Combining this with \eqref{eq4}, one obtains that there exists a $\tilde{\phi}\in\sD$ such that
\beg{align}\label{eqD}
&|\na_{b_t(x)}u_t(x)+b_t(x)-\na_{b_t(y)}u_t(y)-b_t(y)|\nonumber\\
&\le(1+\|\na u\|_{T,\8})\phi(|x-y|)+\|b\|_{T,\8}\|\na u_t(x)-\na u_t(y)\|\nonumber\\
&\le(1+\|\na u\|_{T,\8})\phi(|x-y|)+2\|b\|_{T,\8}\sq{|x-y|}I_{\{|x-y|\ge1\}}\nonumber\\
&~~+c\|b\|_{T,\8}\|\na_{b}u+b\|_{T,\8}\int_s^T\e^{-\la(t-s)}\Big(\ff{|x-y|}{(t-s)}\wedge\ff{1}{\sq{t-s}}\Big)I_{\{|x-y|\le1\}}\d t\nonumber\\
&\le(1+\|\na u\|_{T,\8})\phi(|x-y|)+2\|b\|_{T,\8}\sq{|x-y|}I_{\{|x-y|\ge1\}}\nonumber\\\nonumber
&~~+c\|b\|_{T,\8}\|\na_{b}u+b\|_{T,\8}|x-y|\log(\e+\ff{1}{|x-y|})I_{\{|x-y|\le1\}}\nonumber\\
&\le c\sq{\phi^2(|x-y|)+|x-y|}\nonumber\\
&=: \tilde{\phi}(|x-y|), 
\end{align}
the last inequality was due to the fact that for $x\in[0,1],~\sq{x}\log(\e+\ff{1}{x})$ is an increasing function.

Using $\|\na u\|_{T,\8}+\|b\|_{T,\8}<\8$,  \eqref{eq4}, \eqref{eq2.2} and \eqref{eqD}, we derive
\beg{align}\label{eqh}
\|\na^2u\|_{T,\8}&=\int_0^T\e^{-r\la}\sup_{x\in\R^d}\|\na^2 P_{0,r}^0\{\na_{b_r}u_r+b_r\}(x)\|\d r\nonumber\\
&\le\int_0^T\e^{-r\la}\ff{c_1\tilde{\phi}(c_2\sq{r})}{r}\d r
=:\de_{\tilde{\phi}}(\la).
\end{align}
Noting that $\tilde{\phi}\in \sD$, we have $\int_0^T\ff{c_1\tilde{\phi}(c_2\sq{r}\|)}{r}\d r<\8$, which implies that $\de_{\tilde{\phi}}(\la)\rightarrow 0$ as $\la\rightarrow\8$.
\end{proof}
We provide the regularization representation \eqref{eq2} of  solution to \eqref{eq1}.
\beg{lem}\label{lem2}
Assume (A1) and (A2'). Then for any $T>0$, there exists a constant $\la(T)>0$ such that  for any $\la\ge\la(T)$, it  holds that $\P$-a.s.
\beg{align}\label{e2}
X_t&=X_0+u_0(X_0)-u_t(X_t)+\int_0^t\{\si_s+(\na u_s)\si_s\}(X_s)\d W_s\nonumber\\
&~~+\int_0^t\{\la u_s+B_s\}(X_s)\d s,
\end{align}
where $u$ solves \eqref{eq4}.
\end{lem}
\beg{proof}
Let $G_r=\na_{b_r}u_r+b_r$, $r\ge0$. For fixed $\de>0$, let
\beg{align*}
F^{(\de)}_{s,r}(x)=P_{s,r+\de}^0G_r(x),~~0\le s< r\le T, x\in\R^d.
\end{align*}
According to (A2) and \eqref{u1}, we know $G_r$ is bounded and measurable. Then, we obtain from \eqref{ine-na2} that
\beg{align}\label{ineqF}
\sup_{0\le s< r\le T}\{\|F^{(\de)}_{s,r}\|_\8+\|\na F^{(\de)}_{s,r}\|_\8+\|\na^2 F^{(\de)}_{s,r}\|_\8\}<\8.
\end{align}
By \eqref{eq3} and It\^o's formula, we derive that for any $0\le s\le r\le T$
\beg{align*}
\d F^{(\de)}_{s,r}(Z_{r,t}^x)=\tilde{L}_t F^{(\de)}_{s,r}(Z_{r,t}^x)\d t+\<\na F^{(\de)}_{s,r}(Z_{r,t}^x),\si_t(Z_{r,t}^x)\d W_t\>, t\ge r,
\end{align*}
which yields that 
\beg{align}\label{d-F}
\ff{\d}{\d s}F^{(\de)}_{s,r}(x):&=-\lim_{v\da0}\ff{F^{(\de)}_{s-v,r}(x)-F^{(\de)}_{s,r}(x)}{v}
=-\lim_{v\da0}\ff{P_{s-v,s}^0P_{s,r+\de}^0G_r(x)-F^{(\de)}_{s,r}(x)}{v}\nonumber\\
&=-\lim_{v\da0}\ff{\E P_{s,r+\de}^0G_r(Z_{s-v,s}^x)-F^{(\de)}_{s,r}(x)}{v}
=-\lim_{v\da0}\ff{\E F^{(\de)}_{s,r}(Z_{s-v,s}^x)-F^{(\de)}_{s,r}(x)}{v}\nonumber\\
&=-\lim_{v\da0}\ff{1}{v}\E\int_{s-v}^s(\tilde{L}_tF^{(\de)}_{s,r})(Z_{s-v,t}^x)\d t
=-\tilde{L}_sF^{(\de)}_{s,r}(x), r>0, \mbox{a.e.} ~s\in[0,r].
\end{align}
Let 
\beg{align}\label{u-de}
u_s^{(\de)}=\int_s^T\e^{-\la(t-s)}P_{s,t+\de}^0G_t\d t=\int_s^T\e^{-\la(t-s)}F^{(\de)}_{s,t}\d t, s\in[0,T].
\end{align}
Then we obtain from \eqref{u-de}, \eqref{ineqF} and \eqref{d-F} that
\beg{align*} 
\pp_su_s^{(\de)}=(\la-\tilde{L}_s)u_s^{(\de)}-P_{s,s+\de}^0(\na_{b_s}u_s+b_s).
\end{align*}
By It\^o's formula, we arrive at 
\beg{align}\label{Ito-u-de}
\d u_s^{(\de)}(X_s)&=\{\tilde{L}_su_s^{(\de)}+\na_{b_s}u_s^{(\de)}+\pp_su_s^{(\de)}\}(X_s)\d s+\<\na u_s^{(\de)}(X_s),\si_s(X_s)\d W_s\>\nonumber\\
&=\left\{ \la u_s^{(\de)}+\na_{b_s}u_s^{(\de)}-P_{s,s+\de}^0\na_{b_s}u_s\right\}(X_s)\d s\nonumber\\
&\qquad +\<\na u_s^{(\de)}(X_s),\si_s(X_s)\d W_s\>.
\end{align}

It follows from \eqref{u-de}  and \eqref{eq4}  that
\beg{align}\label{u-de-u}
u_s^{(\de)}-u_s=\int_s^T\e^{-\la(t-s)}(P_{s,t}^0\{P_{t,t+\de}^0G_t-G_t\})\d t,
s\in[0,T].
\end{align}
By \eqref{eqD}, $G_t(\cdot)$ is continuous. Then
$$\lim_{\de\ra 0^+}P_{s,t+\de}^0G_t=P^0_{s,t}G_t,$$ 
which, together with the boundedness of $\|\na u\|$ and $b$, implies by  the dominated convergence theorem that
\beg{align}\label{u-de-u-1}
\lim_{\de\ra 0^+}\left|u_s^{(\de)}-u_s\right|&\leq \int_s^T\e^{-\la(t-s)}\lim_{\de\ra 0^+}\left|P_{s,t}^0\{P_{t,t+\de}^0G_t-G_t\})\right|\d t\nonumber\\
&=0,~s\in[0,T].
\end{align}
By using the boundedness of $\|\na u\|$ and $b$ again, we can derive from \eqref{eq10} and \eqref{u-de} that $\sup_{\de\in (0,1)}\|\na u^{(\de)}\|_{T,\infty}<\infty$. Moreover, combining \eqref{u-de-u} with \eqref{eq10}, we obtain from the dominated convergence theorem that
\beg{align}\label{de-u}
\lim_{\de\da0}\|\na u_s^{(\de)}-\na u_s\|
&=\lim_{\de\da0}\|\int_s^T\e^{-\la(t-s)}\na P_{s,t}^0\{P_{t,t+\de}^0G_t-G_t\}\d t\|\nonumber\\
&\le\lim_{\de\da0}\|\int_s^T\ff{c\e^{-\la(t-s)}}{\sq{(t-s)}}\sq{P_{s,t}^0|P_{t,t+\de}^0G_t-G_t|^2}\d t\|=0.
\end{align}
Combining this with \eqref{Ito-u-de}, \eqref{u-de-u-1} and \eqref{eq1}, we obtain \eqref{e2}.
\end{proof}

{\bf Proof of Theorem \ref{Dthm}}
\beg{proof}
(1) By Lemma \ref{lu}, we can take $\la(T)>0$ large enough such that for any $\la\ge\la(T)$, the unique solution $u$ to \eqref{eq4} satisfies 
 \beg{align}\label{gru}
 \|\na u\|_{T,\8}<\ff{1}{2}.
 \end{align}
This implies that  $\Ph_t(x):=x+u_t(x)$ is a differeomorphism and satisfies that for $(t,x)\in[0,T]\times\R^d$,
\beg{align}\label{grp}
&\ff 1 2\le\|\na\Ph_t(x)\|\le\ff 3 2, ~~
\ff 2 3\le\|\na\Ph_t^{-1}(x)\|\le 2.
\end{align}
Since $u\in C([0,T];C_b^1(\R^d;\R^d))$,  we define $\Phi :C([0,T];\R^d)\rightarrow C([0,T];\R^d)$ as
\beg{align}\label{Ph}
\Ph(\xi)(t)=\Ph_t(\xi_t),~~\xi\in C([0,T];\R^d), t\in[0,T].
\end{align}
Moreover, it follows from \eqref{grp} that 
\beg{align*}
&|\Ph_{t+\tri t}^{-1}(\xi_{t+\tri t})-\Ph_{t}^{-1}(\xi_t)|\\
&\le|\Ph_{t+\tri t}^{-1}(\xi_{t+\tri t})-\Ph_{t+\tri t}^{-1}(\xi_{t})|
+|\Ph_{t+\tri t}^{-1}(\xi_{t})-\Ph_{t}^{-1}(\xi_{t})|\\
&\le\|\na\Ph_{t+\tri t}^{-1}(\cdot)\|_\8|\xi_{t+\tri t}-\xi_t|
+|\Ph_{t+\tri t}^{-1}(\xi_{t+\tri t})-\Ph_{t+\tri t}^{-1}(\Ph_{t+\tri t}(\Ph_{t}^{-1}(\xi_t)))|\\
&\le 2\Big\{|\xi_{t+\tri t}-\xi_t|+|\xi_{t+\tri t}-\Ph_{t+\tri t}(\Ph_{t}^{-1}(\xi_t))|\Big\},
\end{align*}
which yields that $\Ph_\cdot^{-1}(\xi_\cdot)$ is also continuous. Hence $\Ph$ is a homeomorphisms on $C([0,T];\R^d)$ with
\beg{align}\label{iPh}
\Ph^{-1}(\xi)(t)=\Ph_t^{-1}(\xi_t),~~\xi\in C([0,T];\R^d), t\in[0,T].
\end{align}
Then $\Ph$ induces a homeomorphism on $C([0,T];\R^d)\times C([0,T];\R^d)$ defined as in Section 2 (setting $E=C([0,T];\R^d)$) which is still denoted by $\Ph$, and its inverse is still denoted by $\Ph^{-1}$.  Furthermore, it follows from \eqref{grp} and \eqref{Ph} that for any $\xi,\et\in C([0,T];\R^d)$
\beg{align}\label{dis}
&\ff 1 2\rho_T(\xi,\et)\le\rho_T\circ\Ph(\xi,\et)\le\ff 3 2\rho_T(\xi,\et).
\end{align}
These mean that condition \eqref{inq-1} holds for $\rho_T\circ{\Ph}$ by setting $c_1=\ff 1 2,c_2=\ff 3 2 $.

By setting $Y_t=\Ph_t(X_t)$, it follows from Lemma \ref{lem2} that 
\beg{align}\label{equ-YY}
Y_t&= Y_0+\int_0^t (\la u_s+B_s)\circ \Ph_s^{-1}(Y_s)\d s\nonumber\\
&+\int_0^t (\na\Ph_s\si_s)\circ \Ph_s^{-1}(Y_s)\d W_s,~t\in [0,T].
\end{align}
Moreover, it follows from \eqref{nn-b-B} and \eqref{eqh} that
\beg{align}\label{lip-y}
\|\na(\la u_s+B_s)\|_{T,\infty}+\|\na (\na\Ph_s\si_s) \|_{T,\infty}<\infty.
\end{align}
Then there exists a constant $C>0$ (see e.g. \cite[Theorem 1]{U} or \cite{BWY}) such that
\beg{align*}
\W_{2}^{\rho_T}(\Q\circ\Ph^{-1},\P^x\circ\Ph^{-1})^2\le CH(\Q\circ\Ph^{-1}|\P^x\circ\Ph^{-1}).
\end{align*}
Combining  this with Lemma \ref{gthm} and \eqref{dis}, we have that
\beg{align*}
\W_{2}^{\rho_T}(\Q,\P^x)^2\le 2CH(\Q |\P^x ).
\end{align*}

(2)  Based on \cite[Theorem 2.1]{WZ}, it suffices to verify the following assertions respectively: 
\beg{align}\label{q1}
\W_{2}^{\rho_T}(\Q,\P^x)^2&\le c_1H(\Q|\P^x),~~\Q\in\sP(C([0,T];\R^d)),\\
\W_{2}^{\rho_T}(\P^x,\P^y)^2&\le c_2|x-y|^2,~~x,y\in \R^d,\label{q2}
\end{align}
for some     constants $c_1$ and $c_2$.

Since \eqref{q1} has been proved in (1), we only need to prove \eqref{q2}.   Noting that the law of $(X_t^x,X_t^y)_{t\in[0,T]}$ is a coupling of $\P^x$ and $\P^y$, we obtain that
\beg{align*}
\W_{2}^{\rho_T}(\P^x,\P^y)^2\le\E[\rho_T(X^x,X^y)^2]=\E\Big(\sup_{t\in[0,T]}|X_t^x-X_t^y|^2\Big).
\end{align*}
Denote by $Y_t^{\Ph_0(x)}$ the solution of \eqref{equ-YY} with $Y_0=\Ph_0(x)$. By \eqref{lip-y}, it is easy to derive from the B-D-G inequality that
\beg{align*}
\E\sup_{t\in [0,T]}|Y_t^{\Ph_0(x)}-Y_t^{\Ph_0(y)}|^2\leq C|\Ph_0(x)-\Ph_0(y)|^2.
\end{align*}
Combining this with \eqref{dis}, we have that
\beg{align*}
\E\sup_{t\in [0,T]}|X_t^x-X_t^{y}|^2\leq 4\E\sup_{t\in [0,T]}|Y_t^{\Ph_0(x)}-Y_t^{\Ph_0(y)}|^2\leq 9C|x-y|^2.
\end{align*}

\end{proof}

\section{Proof of Theorem \ref{TCI-1}}
For the reader's  convenience, we sketch the construction of homeomorphism $\Phi$. 
To this end, we consider the following elliptic equation
\beg{align}\label{e-b}
(\sL_1^{b_1}+\sL_2^{\si}-\la)u=b_1,
\end{align}
where $\sL_1^{b_1}=\na_{b_1}, ~~\sL_2^{\si}:=\ff{1}{2}\sum_{i,j}\<\si\si^*e_i,e_j\>\na_{e_i}\na_{e_j}$.

Before moving on, we introduce some spaces and  notations. For $(p,\al)\in[1,\8]\times(0,2]-\{\8\}\times\{1\}$, let $H_p^\al:=(I-\De)^{-\ff{\al}{2}}(L^p(\R^d))$ be the usual Bessel potential space with the norm 
\beg{align*}
\|f\|_{\al,p}:=\|(I-\De)^{\ff{\al}{2}}f\|_p\asymp\|f\|_p+\|\De^{\ff{\al}{2}}f\|_p,
\end{align*}
where $\|\cdot\|_p$ is the usual $L^p$-norm in $\R^d$, and $(I-\De)^{\ff{\al}{2}}f$ and $\De^{\ff{\al}{2}}f$ are defined through the Fourier transformation 
$$(I-\De)^{\ff{\al}{2}}f:=\cF^{-1}((1+|\cdot|^2)^{\ff{\al}{2}}\cF f),
\De^{\ff{\al}{2}}f:=\cF^{-1}(|\cdot|^\al\cF f).$$
For $(p,\al)=(\8,1)$, we define $H_\8^1$ as the space of Lipschitz functions with finite norm 
$$\|f\|_{1,\8}:=\|f\|_\8+\|\na f\|_\8<\8.$$
Notice that for $n=1,2$ and $p\in(1,\8)$, an equivalent norm in $H_p^n$ is given by 
$$\|f\|_{n,p}:=\|f\|_p+\|\na^n f\|_p<\8.$$
The above facts are standards and can be found in \cite{T1}.

The following Lemma shows the solvability of   equation \eqref{e-b}, which is a consequence of  \cite[Theorem 7.5]{XZ}.
\beg{lem}\label{s-u}
Suppose that $(H^\si)$ holds and $b_1\in L^p(\R^d)$ for some $p>d$. Then there exists sufficient large constant $\la_1$ such that  for all $\la\ge\la_1$ there exists a unique solution $u\in H_p^2$ to equation \eqref{e-b}. Moreover,  for $v\in(0,2)$ with $\ff{d}{p}<2-v$, we have
\beg{align}\label{u-L}
\la^{\ff{1}{2}(2-v)}\|u\|_{v,p}+\|\na^2u\|_p\le c\|b_1\|_p.
\end{align}
\end{lem}

Recall the following Sobolev embedding for $p\in[1,\8], \al\in[0,2]$
\beg{align*}
H_p^\al\subset H_\8^{\al-\ff{d}{p}}\subset C_b^{\al-\ff{d}{p}},~~~\al p>d.
\end{align*}
Combining this with \eqref{u-L}, one can see  that there exist $c,\la_1\ge1$ such that for all $\la\ge\la_1$,
\beg{align}\label{u-u-L}
\|u\|_\8+\|\na u\|_\8\le c\la^{\ff{1}{2}(\ff{d}{p}-1)}.
\end{align}
Define $\Phi(x):=x+u(x)$. By \eqref{u-u-L} with $\la$ large enough, the map $x\ra\Phi(x)$ forms a $C^1$-diffeomorphism and 
\beg{align}\label{phi-up}
\ff{1}{2}\le \|\na\Phi\|_\8, \|\na\Phi^{-1}\|_\8\le 2.
\end{align}
The following Lemma present the regular representation of solution to \eqref{b-S} by Zvonkin's transformation. This result is due to \cite[Lemma 7.6]{XZ}.
\beg{lem}\label{x-phi}
$X_t$ solves SDE \eqref{b-S} if and only if $Y_t:=\Phi(X_t)$ solves 
\beg{align}\label{phi-r}
\d Y_t=\til{b}(Y_t)\d t+\til{\si}(Y_t)\d W_t,
\end{align}
with initial value $y:=\Phi(x)$  and 
\beg{align*}
\til{b}(y):=(\la u+\na\Phi\cdot b_2)\circ\Phi^{-1}(y), ~~\til{\si}(y):=(\na\Phi\cdot \si)\circ\Phi^{-1}(y).
\end{align*}
\end{lem}
The following Lemma shows that the conditions for $b_2$ in $(H^b)$ are preserved under Zvomnkin's transformation.
\beg{lem}\label{H-B}
Under $(H^b)$, Then for $\la$ large enough, 
\beg{enumerate}
\item[($\tld H^b_2$)] $\tld b$ satisfies one of the following conditions\\
(1) there exist $r>-1$, $\tld\ka_1>0$, $\tld\ka_2\geq 0$ and $\tld\ka_3\geq 0$ such that 
\beg{align}\label{b-dis}
\<\tld b(y),y\>\leq -\tld \ka_1|y|^{2+r}+\tld \ka_2,\qquad |\tld b(y)|\leq \tld \ka_3(1+|y|^{r+1}),~y\in\R^d;
\end{align}
(2) there exist $\ka_4\geq 0$ such that
\beg{align}\label{b-lin}
|\tld b(y)|\leq \tld\ka_4(1+|y|),~y\in\R^d.
\end{align}
\end{enumerate}
\end{lem}

We establish $T_1(C)$  by ``Gaussian tail'' following \cite[Theorem 2.3]{DGW}, and we recall the following lemma there.
\beg{lem}\label{G}
A given probability measure $\mu$ on $(E,\rho)$ satisfies the $L^1$-transportation cost information inequality with some constant $C$ if and only if 
\beg{align}\label{G-tail}
\int\int\e^{\de \rho^2(x,y)}\d\mu(x)\d\mu(y)<+\8, ~~\de\in(0,\ff{1}{4C}),
\end{align}
holds.
\end{lem}

The following two lemmas contribute to establishing \eqref{G-tail} for solutions $Y_t$ to the equation \eqref{phi-r}. 
By the definition of $\tld \si$ and \eqref{phi-up}, it is clear that
\beg{align}\label{si-hs}
\sup_{y\in \R^d} \|\tld\si(y)\|_{HS}<\infty.
\end{align} 

\beg{lem}\label{exp-int}
Assume that $\tld b$ satisfies ($\tld H^b_2$). Then 
$$\E\exp\left\{\la\int_0^T|Y_t|^{2r+2}\d t\right\}<\infty,$$
where, under the condition (1) in ($\tld H^b_2$)
$$\la <  2^{-(r-1)^-}  \tld \ka_1^2  \|\tld \si\|_\infty^{-2}, ~\text{and}~ \la \leq  2^{-(r-1)^-}  \tld \ka_1^2  \|\tld \si\|_\infty^{-2}~\text{even if}~r<0;$$
under the condition (2) in ($\tld H^b_2$)
$$\la\leq \ff {e^{-(2+3\tld\ka_4T)}} {2\|\tld\si\|^2_\infty}.$$
\end{lem}
\beg{proof}
We first prove this lemma under the condition (1). It follows from  It\^o's formula that
\beg{align}\label{Ito-Y1}
\d |Y_t|^2\leq \left(-2\tld\ka_1|Y_t|^{r+2}+2\tld \ka_2\right)\d t+2\<Y_t,\tld \si(Y_t)\d W_t\>+\|\tld \si(Y_t)\|^2_{HS}\d t.
\end{align}

For $r\in (-1,0)$, we have by \eqref{Ito-Y1} that
\beg{align}\label{Ito-Yr2}
\d \left(1+|Y_t|^2\right)^{\ff {r+2} 2}&\leq -\ff {r+2} 2\left(1+|Y_t|^2\right)^{\ff r 2}\left(-2\tld \ka_1|Y_t|^{2+r}+2\tld \ka_2+\|\tld\si\|^2_{HS}\right)\d t\nonumber\\
&\qquad +(r+2)\left(1+|Y_t|^2\right)^{\ff r 2}\<Y_t,\tld\si(Y_t)\d W_t\>\nonumber\\
&\qquad +\ff {r(r+2)} 2 \left(1+|Y_t|^2\right)^{\ff r 2-1}|\tld\si^*(Y_t)Y_t|^2\d t \nonumber\\
&\leq -2^{\ff r 2}(r+2)\tld \ka_1|Y_t|^{2+2r}\d t+(r+2\left(2^{\ff r 2}+\tld \ka_1+\|\tld\si\|^2_{HS,\infty}\right)\d t\nonumber\\
&\qquad +(r+2)\left(1+|Y_t|^2\right)^{\ff r 2}\<Y_t,\tld\si(Y_t)\d W_t\>,
\end{align} 
where in the last inequality we use
$$(1+y^2)^{\ff r 2}y^{2+r}\geq 2^{\ff r 2}y^{2+2r}-2^{\ff r 2},~y\geq 0,~r\in(-1,0).$$
Let 
$$\ta_n=\inf\{ t>0~|~\int_0^t |Y_s|^{2+2r}\d s\}.$$
Then it follows from \eqref{Ito-Yr2} that 
\beg{align*}
&\E \exp\left\{\la\int_0^{T\we\ta_n} |Y_t|^{2+2r}\d t\right\}\\
&\leq C_\la\E\exp\left\{\ff {\la} {\tld \ka_1 2^{\ff r 2}}\int_0^{T\we\ta_n}(1+|Y_t|^2)^{\ff r 2}\<Y_t,\tld\si(Y_t)\d W_t\>\right\}\\
&\leq C_\la\left(\E\exp\left\{\ff {2^{1-\ff r 2}\la} {\tld \ka_1 }M_{T\we\ta_n}-\ff {2^{1-r}\la^2} {\tld \ka_1^2 }\<M\>_{T\we\ta_n}\right\}\right)^{\ff 1 2}\left(\E\left\{\ff {2^{1-r}\la^2 } {\tld \ka_1^2}\<M\>_{T\we\ta_n}\right\}\right)^{\ff 1 2}\\
&\leq C_\la\left(\E\left\{\ff {2^{1-r}\la^2\|\tld\si\|^2_\infty} {\tld \ka_1^2}\int_0^{T\we\ta_n} |Y_t|^{2+2r}\d t\right\}\right)^{\ff 1 2},
\end{align*}
where
\beg{align*}
C_\la & =\exp\left\{-\ff {\la(1+|Y_0|^2)^{\ff {r+2} 2}+\la(r+2)\left(2^{\ff r 2-1}\tld \ka_1+\tld \ka_2+\|\tld\si\|^2_{HS,\infty}\right)T} {\tld \ka_1(r+2)2^{\ff r 2}}\right\}\\
M_{T\we\ta_n}& =\int_0^{T\we\ta_n}(1+|Y_t|^2)^{\ff r 2}\<Y_t,\tld\si(Y_t)\d W_t\>.
\end{align*}
Choosing $\la=\ff {\tld \ka_1^2} {2^{1-r}\|\tld \si\|_\infty^2}$ and letting $n\ra +\infty$, we have
\beg{align*}
\E \exp\left\{\ff {\tld \ka_1^2} {2^{1-r}\|\tld \si\|_\infty^2}\int_0^{T} |Y_t|^{2+2r}\d t\right\}<\infty.
\end{align*}

Let $\ep_1\in (0,1)$ if $r>0$ and $\ep_1=0$ if $r=0$. For $r\geq 0$, it follows from It\^o's formula and the H\"older inequality that
\beg{align*}
\d |Y_t|^{r+2}& \leq -(r+2)\tld\ka_1 |Y_t|^{2r+2}\d t+\ff {r+2} 2\left(2\tld\ka_2+\|\tld\si\|^2_{HS,\infty}+r\|\tld\si\|^2_{\infty}\right)|Y_t|^{r}\d t\\
&\qquad + (r+2)|Y_t|^{r}\<Y_t,\tld\si(Y_t)\d W_t\>\\
&\leq -(r+2)\tld\ka_1 |Y_t|^{2r+2}\d t+\ep_1(r+2)\tld\ka_1|Y_t|^{2r+2}\d t\\
&\qquad +\ff {(r+2)^2\left(2\tld \ka_2+(r+1)\|\tld \si\|^2_{HS,\infty}\right)\left(2r\tld\ka_2+r(r+1)\|\tld\si\|^2_{HS,\infty}\right)^{\ff r {r+2}}} {4(r+1)(4\ep_1\tld \ka_1(r+1))^{\ff r {r+2}}} \d t\\
&\qquad +(r+2)|Y_t|^r\<Y_t,\tld\si(Y_t)\d W_t\>.
\end{align*}
Then for  $\la=\ff {(1-\ep_1)^2\tld\ka_1^2} {2\|\tld\si\|^2_\infty}$, there is $C_\la>0$ such that
\beg{align*}
\E\exp\left\{\la\int_0^T|Y_t|^{2r+2}\d t\right\}&\leq C_\la\E\exp\left\{\ff {\la} {(1-\ep_1)\tld\ka_1}\int_0^T|Y_t|^r\<Y_t,\tld\si(Y_t)\d W_t\>\right\}\\
&\leq  C_\la\left(\E\exp\left\{\ff {2\la^2\|\tld\si\|^2_\infty} {(1-\ep_1)^2\tld\ka_1^2}\int_0^T|Y_t|^{2r+2}\d t \right\}\right)^{\ff 1 2}\\
&=C_\la \left(\E\exp\left\{\la\int_0^T|Y_t|^{2r+2}\d t \right\}\right)^{\ff 1 2}
\end{align*}
which yields that for any $\la<\ff 1 2  \tld\ka_1^2 \|\tld\si\|^{-2}_\infty$ if $r>0$ and $\la\leq \ff 1 2  \tld\ka_1^2 \|\tld\si\|^{-2}_\infty$  if $r=0$, 
$$\E\exp\left\{\la\int_0^T|Y_t|^{2r+2}\d t\right\}<\infty.$$

Next, we prove this lemma under the condition (2). It follows from It\^o's formula and the H\"older inequality  that for $\al>3\tld\ka_4$
\beg{align*}
\d \left(e^{-\al t}|Y_t|^2\right)&\leq e^{-\al t}(\tld\ka_4+\|\si\|^2_{HS,\infty})\d t-(\al-3\tld\ka_4)e^{-\al t}|Y_t|^2\d t\\
&\qquad +2e^{-\al t}\<Y_t,\tld\si(Y_t)\d W_t\>.
\end{align*}
Then
\beg{align*}
&\exp\left\{\la \int_0^T e^{-\al t}|Y_t|^2-\ff {\la(\tld \ka_4+\|\tld\si\|^2_{HS,\infty})(1-e^{-\al T})} {\la(\al-3\tld\ka_4)}\right\}\\
&\qquad \leq \exp\left\{\ff {2\la} {\al-3\tld\ka_4}\int_0^Te^{-\al t}\<Y_t,\tld\si(Y_t)\d W_t\>\right\}.
\end{align*}
Let $\la=\ff {(\al-3\tld\ka_4)^2} {8\|\tld\si\|^2_\infty}$. Then we have that
$$\E\exp\left\{\ff {(\al-3\tld\ka_4)^2} {8\|\tld\si\|^2_\infty}\int_0^T e^{-\al t}|Y_t|^2\d t\right\}<\infty.$$
By choosing the optimal $\al=\ff 2 T +3\tld\ka_4$, we have 
\beg{align*}
\E\exp\left\{\ff {e^{-(2+3\tld\ka_4T)}} {2\|\tld\si\|^2_\infty}\int_0^T  |Y_t|^2\d t\right\}<\infty.
\end{align*}

\end{proof}

Let $W^{(1)}_t, {W^{(2)}_t}$ be  two independent Brownian motions defined on the filtered probability $(\Om,\sF,(\sF_t),\P)$, and $Y^{(i)}_t$, $i=1,2$ are  solutions of \eqref{phi-r} driven by $W^{(i)}_t$ with the same initial value $y$.  Let
$$Z_t=Y^{(1)}_t-Y^{(2)}_t.$$
Then
$$\E\exp\left\{\de \sup_{t\in [0,T]}|Z_t|^2\right\}=\iint_{C([0,T];\R^d)\times C([0,T];\R^d)}\e^{\de\rh_T(\xi,\et)}\tld\P^{y}(\d \xi)\tld\P^{y}(\d\et),$$
where $\tld\P^{y}$ is the law of $Y^{(1)}$ on $C([0,T];\R^d)$. 
\beg{lem}\label{exp-int2}
Suppose the assumptions in Lemma \ref{exp-int} hold. Let 
$$\de<\ff  1 {4\|\tld\si\|^2_\infty T\left(2+\tld\ka_3^2 \tld\ka_1^{-2}2^{(r-1)^-}  \right)}$$
if $\tld b$ satisfies \eqref{b-dis}, and let 
$$\de<\ff  1 {8 \|\tld\si\|^2_\infty T\left(1+\tld\ka_4^2 \e^{2+3\tld\ka_4T}\right)}$$
if $\tld b$ satisfies \eqref{b-lin}. 
Then
$$\E\exp\left\{\de \sup_{t\in [0,T]}|Z_t|^2\right\}<\infty.$$
\end{lem}
\beg{proof}
It follows from It\^o's formula that
\beg{align}\label{Ito-Z}
\d \sq{1+|Z_t|^2}&=\ff {\<\tld b(Y^{(1)}_t)-\tld b(Y^{(2)}_t),Y^{(1)}_t-Y^{(2)}_t\>} {\sq{1+|Z_t|^2}}\d t+ \ff {\|\tld\si(Y^{(1)}_t)\|^2_{HS}+\|\tld\si(Y^{(2)}_t)\|^2_{HS}} {\sq{1+|Z_t|^2}}\d t\nonumber\\
&\qquad +\ff {\<Z_t,\tld\si(Y^{(1)}_t)\d W^{(1)}_t-\tld\si(Y^{(2)}_t)\d W^{(2)}_t\>} {\sq{1+|Z_t|^2}}\nonumber\\
&\qquad -\ff {|\tld\si^*(Y^{(1)}_t)Z_t|^2+|\tld\si^*(Y^{(2)}_t)Z_t|^2} {\left(1+|Z_t|^2\right)^{\ff 3 2}}\d t.
\end{align}

If $\tld b$ satisfies \eqref{b-lin}, then
\beg{align*}
\ff {\<\tld b(Y^{(1)}_t)-\tld b(Y^{(2)}_t),Y^{(1)}_t-Y^{(2)}_t\>} {\sq{1+|Z_t|^2}}\leq \tld\ka_4\left(2+|Y^{(1)}_t|+|Y^{(2)}_t|\right).
\end{align*}
Putting this into \eqref{Ito-Z}, we have that
\beg{align*}
\d \sq{1+|Z_t|^2}&\leq 2(\tld \ka_4+\|\tld\si\|^2_{HS,\infty})\d t+ \tld\ka_4\left(|Y^{(1)}_t|+|Y^{(2)}_t|\right)\d t+\d M_t
\end{align*}
with
$$M_t=\int_0^t\ff {\<Z_s,\tld\si(Y^{(1)}_s)\d W^{(1)}_s-\tld\si(Y^{(2)}_s)\d W^{(2)}_s\>} {\sq{1+|Z_t|^2}}.$$
Then 
\beg{align*}
&\E\exp\left\{\be\left(\sup_{t\in [0,T]}\sq{1+|Z_t|^2}-2(\tld\ka_4+\|\tld\si\|^2_{HS,\infty})\right)\right\}\\
&\qquad \leq \E\left( \sup_{t\in[0,T]}\e^{\be M_t}\e^{\be\tld\ka_4\int_0^T\left(|Y^{(1)}_t|+|Y^{(2)}_t|\right)\d t}\right)\\
&\qquad \leq \left(\E\sup_{t\in [0,T]}\e^{2\be M_t}\right)^{\ff 1 2}\left(\E\e^{2\be\tld\ka_4\int_0^T\left(|Y^{(1)}_t|+|Y^{(2)}_t|\right)\d t}\right)^{\ff 1 2}\\
&\qquad \leq 2\left(\E \e^{2\be M_T}\right)^{\ff 1 2}\left(\E\e^{2\be\tld\ka_4\int_0^T |Y^{(1)}_t| \d t}\right) 
\end{align*}
where  we use the Doob's maximal inequality and the independence of $Y^{(1)}_t$ and $Y^{(2)}_t$. Since
\beg{align*}
\<M\>_T=\int_0^T\ff {|\tld\si^*(Y^{(1)}_t)Z_t|^2+|\tld\si^*(Y^{(2)}_t)Z_t|^2} {1+|Z_t|^2}\d t\leq 2\|\tld\si\|^2_\infty T,
\end{align*}
we have
\beg{align}\label{exp-m-1}
\E \e^{2\be M_T}\leq \left(\E\e^{2\be M_T-2\be^2\<M\>_T}\right)\e^{4\be^2\|\tld\si\|^2_\infty T}=\e^{4\be^2\|\tld\si\|^2_\infty T}.
\end{align}
By the H\"older inequality, we have that
\beg{align*}
\E\e^{2\be\tld\ka_4\int_0^T |Y^{(1)}_t| \d t}\leq \e^{2\be^2\tld\ka_4^2\|\tld\si\|^2_\infty T\e^{2+3\tld\ka_4T}}\E\e^{\ff {e^{-(2+3\tld\ka_4T)}} {2\|\tld\si\|^2_\infty}\int_0^T |Y^{(1)}_t|^2 \d t}.
\end{align*}
Combining this with \eqref{exp-m-1}, we can derive that
\beg{align*}
&\E\exp\left\{\be\left(\sup_{t\in [0,T]}\sq{1+|Z_t|^2}-2(\tld\ka_4+\|\tld\si\|^2_{HS,\infty})-1\right)\right\}\\
&\qquad \leq 2\exp\left\{2\be^2\|\tld\si\|^2_\infty T\left(1+\tld\ka_4^2 \e^{2+3\tld\ka_4T}\right)\right\}\left(\E\e^{\ff {e^{-(2+3\tld\ka_4T)}} {2\|\tld\si\|^2_\infty}\int_0^T |Y^{(1)}_t|^2 \d t}\right).
\end{align*}
Then by Chebychev’s inequality and an optimization of $\be$
\beg{align*}
&\P\left\{\sup_{t\in [0,T]}\sq{1+|Z_t|^2}\geq 1+2(\tld\ka_4+\|\tld\si\|^2_{HS,\infty})+x\right\}\\
&\qquad  \leq 2\left(\E\e^{\ff {e^{-(2+3\tld\ka_4T)}} {2\|\tld\si\|^2_\infty}\int_0^T |Y^{(1)}_t|^2 \d t}\right)\exp\left\{-\ff {x^2} {8\|\tld\si\|^2_\infty T\left(1+\tld\ka_4^2 \e^{2+3\tld\ka_4T}\right)}\right\},
\end{align*}
which yields that 
\beg{align*}
\E\exp\left\{\de \sup_{t\in [0,T]}|Z_t|^2\right\}<\infty,~\de<\ff  1 {8 \|\tld\si\|^2_\infty T\left(1+\tld\ka_4^2 \e^{2+3\tld\ka_4T}\right)}.
\end{align*}

If $\tld b$ satisfies \eqref{b-dis}, then
\beg{align*}
\ff {\<\tld b(Y^{(1)}_t)-\tld b(Y^{(2)}_t),Y^{(1)}_t-Y^{(2)}_t\>} {\sq{1+|Z_t|^2}}\leq \tld\ka_3\left(2+|Y^{(1)}_t|^{r+1}+|Y^{(2)}_t|^{r+1}\right).
\end{align*}
Putting this into \eqref{Ito-Z}, we have that
\beg{align*}
\d \sq{1+|Z_t|^2}&\leq 2(\tld \ka_3+\|\tld\si\|^2_{HS,\infty})\d t+ \tld\ka_3\left(|Y^{(1)}_t|^{r+1}+|Y^{(2)}_t|^{r+1}\right)\d t+\d M_t.
\end{align*}
Since for any $\be>0$ and $0<\tld\be<\ff {  \tld \ka_1^2} {2^{(r-1)^-}  \|\tld \si\|_\infty^{2}}$, we have that
\beg{align*}
\E\e^{2\be\tld\ka_3\int_0^T |Y^{(1)}_t|^{r+1} \d t}\leq \e^{\ff {\be^2\tld\ka_3^2T} {\tld\be }}\E\e^{\tld\be\int_0^T |Y^{(1)}_t|^{2r+2} \d t}<\infty.
\end{align*}
Then
\beg{align*}
&\E\exp\left\{\be\left(\sup_{t\in [0,T]}\sq{1+|Z_t|^2}-2(\tld\ka_3+\|\tld\si\|^2_{HS,\infty})-1\right)\right\}\\
&\qquad \leq 2\exp\left\{2\be^2\|\tld\si\|^2_\infty T+\ff {\be^2\tld\ka_3^2T} {\tld\be }\right\}\left(\E\e^{\tld\be\int_0^T |Y^{(1)}_t|^{2r+2} \d t}\right)
\end{align*}
and
\beg{align*}
&\P\left\{\sup_{t\in [0,T]}\sq{1+|Z_t|^2}\geq 1+2(\tld\ka_3+\|\tld\si\|^2_{HS,\infty})+x\right\}\\
&\qquad  \leq 2\left(\E\e^{\tld\be\int_0^T |Y^{(1)}_t|^{2r+2} \d t}\right)\exp\left\{-\ff {x^2} {8\|\tld\si\|^2_\infty T+4\tld\ka_3^2\tld\be^{-1}T}\right\}.
\end{align*}
Since $\tld\be$ is arbitrary in $(0,\ff {  \tld \ka_1^2} {2^{(r-1)^-}  \|\tld \si\|_\infty^{2}})$, we have that
\beg{align*}
\E\exp\left\{\de \sup_{t\in [0,T]}|Z_t|^2\right\}<\infty,~\de<\ff  1 {4\|\tld\si\|^2_\infty T\left(2+\tld\ka_3^2 \tld\ka_1^{-2}2^{(r-1)^-}  \right)}.
\end{align*}

\end{proof}

{{\bf Proof of Theorem \ref{TCI-1}}
\beg{proof}
Taking the similar arguments as in the proof of Theorem \ref{Dthm}, the assertions of this theorem follows from Lemma  \ref{exp-int2}, Lemma \ref{G} and \eqref{phi-up}. It follows from \eqref{phi-up} that $\Ph$ induces a homeomorphism on $C([0,T];\R^d)$ by using the same argument in Theorem \ref{Dthm}. Moreover, 
\beg{align}\label{distance}
\ff{1}{2}\rho_T(\xi,\et)\le\rho_T\circ\Ph(\xi,\et)\le 2\rho_T(\xi,\et), ~~\xi,\et\in C([0,T];\R^d).
\end{align}
Since $Y_t=\Ph_t(X_t)$, the law of $Y$ is $\P^x\circ \Ph^{-1}$. Then by Lemma \ref{exp-int2} and Lemma \ref{G}, there is a constant $C>0$ such that for any  measure  $Q$  on $C([0,T];\R^d)$, 
\beg{align*}
\W_{1}^{\rho_T}(\Q\circ\Ph^{-1},\P^x\circ\Ph^{-1})\le \sq{CH(\Q\circ\Ph^{-1}|\P^x\circ\Ph^{-1})}.
\end{align*}
Then, following from \eqref{distance} and Lemma \ref{gthm}, 
\beg{align*}
\W_{1}^{\rho_T}(\Q,\P^x)\le \sq{4CH(\Q |\P^x)}.
\end{align*}

\end{proof}

\noindent\textbf{Acknowledgements}

\medskip

The third author was supported by the disciplinary funding of Central University of Finance and Economics, and the National Natural Science Foundation of China (Grant No. 11901604, 11771326).

\end{document}